\documentclass[11pt]{article}
\usepackage{amsfonts}
\usepackage{amsmath}
\usepackage{epsfig}

\newcommand{\C}{\mathbb{C}}

\newcommand{\R}{\mathbb{R}}
\newcommand{\Q}{\mathbb{Q}}
\newcommand{\Z}{\mathbb{Z}}
\newcommand{\PP}{\mathbb{P}}
\newcommand{\TT}{\mathbb{T}}

\newtheorem{thm}{Theorem}
\newtheorem{defn}[thm]{Definition}

\newtheorem{conj}[thm]{Conjecture}

\newtheorem{obs}[thm]{Observation}

\newcommand{\dd}{\;\mathrm{d}}
\newcommand{\e}{\mathrm{e}}
\newcommand{\res}{\mathrm{Res}}

\newcommand{\Li}{\mathrm{Li}} 
\newcommand{\re}{\mathop{\mathrm{Re}}} 
\newcommand{\im}{\mathop{\mathrm{Im}}} 
\newcommand{\ii}{\mathrm{i}} 
\newcommand{\Lf}{\mathrm{L}} 

\newenvironment{proof}{\noindent {\bf Proof.}}{$\Box$}
\newenvironment{acknowledgements}{\noindent{\bf Acknowledgments}\bigskip}{}
\newenvironment{keyword}{\noindent{\bf Keywords}$\,$   }{}
\newenvironment{classification}{\noindent{\bf Classification $\,$   }}{}

\begin{document}

\title{An algebraic integration for Mahler measure}

\author{Matilde N. Lal\'{\i}n\footnote{E-mail address: \texttt{mlalin@math.ubc.ca}}\\ {\small \it  Institut des Hautes \'Etudes Scientifiques } \\{ \small \it Le Bois-Marie, 35, route de Chartres, F-91440 Bures-sur-Yvette, France}}
\date{}

\maketitle

\vspace{-1cm}

\begin{abstract}
There are many examples of several-variable polynomials whose Mahler measure is expressed in terms of special values of polylogarithms. These examples are expected to be related to computations of regulators, as observed by Deninger, and later Rodriguez-Villegas, and Maillot. While Rodriguez-Villegas made this relationship explicit for the two variable case, it is our goal to understand the three variable case and shed some light on the examples with more variables.
\end{abstract}



\maketitle

\section{Introduction}

The (logarithmic) Mahler measure of a Laurent polynomial $P \in \C[x_1^{\pm1}, \dots, x_n^{\pm1}]$ is defined by
\begin{eqnarray}
m(P)& := & \int_0^1 \dots \int_0^1 \log |P(\e^{2 \pi \ii \theta_1}, \dots , \e^{2 \pi \ii \theta_n}) | \dd \theta_1 \dots \dd \theta_n.
\end{eqnarray}

Because of Jensen's formula, there is a simple expression for the Mahler measure in the one-variable case, as a function on the roots of the polynomial. It is natural then, to wonder what happens with several variables.

The problem of finding explicit closed formulas for Mahler measures of several variable polynomials is  hard. However, several examples have been found, especially for two and three variables. Some formulas have been completely proved and some others have been established numerically and are strongly believed to be true.

A remarkable fact is that in most of these examples the Mahler measure of polynomials with integral coefficients can be expressed in terms of special values of L-series or polylogarithms (that is to say, Riemann zeta-functions, Dirichlet L-series, L-series of varieties, zeta functions of number fields, etcetera).

For instance, the first and simplest example in two variables was discovered by Smyth \cite{S1}:
\begin{equation}
m(1+x+y) = \frac{3 \sqrt{3}}{4 \pi} \Lf( \chi_{-3}, 2) = \Lf'(\chi_{-3}, -1)
\end{equation}
where $\chi_{-3}$ is the character of conductor 3.

Another example was computed numerically by Boyd \cite{B2} (and studied by Deninger \cite{D} and Rodriguez-Villegas \cite{RV}),
\begin{equation}
m \left( x + \frac{1}{x} + y + \frac{1}{y} + 1 \right) \stackrel{?}{=} \Lf'(E, 0)
\end{equation}
where $E$ is the elliptic curve of conductor 15 which is the projective closure of the curve $x + \frac{1}{x} + y + \frac{1}{y} + 1=0$ and $\Lf(E, s)$ is the L-function of $E$.

Deninger \cite{D} interpreted computations of Mahler measure in terms of Deligne periods of mixed motives explaining some of the relations to the L-series via Beilinson's conjectures. 

Rodriguez-Villegas \cite{RV} has clarified this relationship by explicitly computing the regulator, and relating this machinery to the cases already (numerically) known by Boyd, proving some of them, and deeply understanding the cases with two variables. Recently Maillot has sketched how one could continue these ideas for more variables. 

It is our goal to develop these ideas and to apply them in order to understand the few known examples with three and more variables involving Dirichlet L-series, Riemann zeta functions and polylogarithms. In this work we describe a general situation and illustrate our explanation with a few examples. In \cite{L3} we will show the computational power of our method by showing how to prove many other formulas of Mahler measures and generalized Mahler measures as well.

\section{Background}
In this section we describe some ingredients that will be used in our construction.

\subsection{Polylogarithms}
The cases that we are going to study involve zeta functions or Dirichlet L-series, but they all may be thought as special values of polylogarithms. In fact, this common feature seems to be the most appropriate way of dealing with the interpretation of these formulas. Here we proceed to recall some definitions and establish some common notation.

\begin{defn} The $n$th polylogarithm is the function defined by the power series
\begin{equation}
 \Li_n(x) := \sum_{k=1}^\infty \frac{x^k}{k^n} \qquad x \in \C, \quad |x| <1.
\end{equation}
\end{defn}
This function can be continued analytically to $\C \setminus (1, \infty)$.
We will work with Zagier's modification of the polylogarithm (\cite{Z2}):
\begin{equation}
 \mathcal{L}_n(x) :=  {\re}_n \left(\sum_{j=0}^{n-1} \frac{2^j B_j}{j!} (\log|x|)^j \Li_{n-j}(x) \right)
\end{equation}
where $B_j$ is the $j$th Bernoulli number and $\re_k$ denotes $\re$ or $\im$ depending on whether $n$ is odd or even. This function is one-valued, continuous in $\PP^1(\C)$, and real analytic in $\PP^1(\C) \setminus \{ 0, 1, \infty\}$.

$\mathcal{L}_n$ satisfies very clean functional equations. The simplest ones are
\[ \mathcal{L}_n\left( \frac{1}{x} \right) = (-1)^{n-1} \mathcal{L}_n(x) \qquad \mathcal{L}_n(\bar{x}) = (-1)^{n-1} \mathcal{L}_n(x).\]

For $n=2$,  one obtains the Bloch Wigner dilogarithm,
\begin{equation}
D(x) = \im(\Li_2(x) -  \log |x| \Li_1(x)) =  \im(\Li_2(x) ) + \log |x| \arg(1-x)
\end{equation}
which satisfies the well-known five-term relation
\begin{equation}
D(x) + D(1-xy) + D(y) + D\left( \frac{1-y}{1-xy} \right) + D\left( \frac{1-x}{1-xy} \right) = 0.
\end{equation}
For $n=3$ we obtain
\begin{equation}
\mathcal{L}_3(x) = \re \left( \Li_3(x) - \log|x| \Li_2(x) + \frac{1}{3} \log^2|x| \Li_1(x) \right).
\end{equation}
This modified trilogarithm satisfies more functional equations, such as the Spence--Kummer relation:
\[\mathcal{L}_3\left(\frac{x(1-y)^2}{y(1-x)^2}\right) + \mathcal{L}_3(xy) + \mathcal{L}_3\left(\frac{x}{y}\right) - 2 \mathcal{L}_3\left(\frac{x(1-y)}{y(1-x)}\right) -2 \mathcal{L}_3\left(\frac{y(1-x)}{y-1}\right)-2\mathcal{L}_3\left(\frac{x(1-y)}{x-1}\right)\]
\begin{equation}\label{SK}
-2 \mathcal{L}_3\left( \frac{1-y}{1-x}\right)-2\mathcal{L}_3(x) -2\mathcal{L}_3(y) + 2 \mathcal{L}_3(1) =0.
\end{equation}



\subsection{Polylogarithmic motivic complexes}\label{Gonchcomplexes}

Given a field $F$, consider $\Z[\PP^1_F]$, the free abelian group generated by the elements of $\PP^1_F$. For each $n$ we are interested in working with this group modulo the (rational) functional equations of the $n$th polylogarithm. Unfortunately, the functional equations of higher polylogarithms are not known explicitly.

For $X$ an algebraic variety, Goncharov \cite{G3,G4,G5}, has constructed some groups that conjecturally correspond to the groups in the above paragraph and they fit into polylogarithmic motivic complexes whose cohomology is related to Bloch groups and is conjectured to be the motivic cohomology of $X$. A regulator can be defined in these complexes and is conjectured to coincide with Beilinson's regulator.

From now on we will follow \cite{G3,G4,G5}. We state definitions and results, the proofs may be found in the mentioned works.

Given a field $F$ one defines inductively some subgroups $\mathcal{R}_n(F)$, then lets
\begin{equation}
{\mathcal B}_n(F) :=\Z[\PP^1_F]/\mathcal{R}_n(F).
\end{equation}
The classes of $x$ in  $\Z[\PP^1_F]$ and in ${\mathcal{B}_n(F)}$ will be denoted by $\{x\}$ and $\{x\}_n$ respectively.
We begin by setting
\begin{equation}
\mathcal{R}_1(F):= \left < \{x\} + \{y\}-\{xy\}; \quad x,y \in F^*, \{0\},\{\infty\}\right >.
\end{equation}
Thus ${\mathcal B}_1(F) = F^*$. Now we proceed to construct a family of morphisms;
\[ \Z[\PP^1_F] \stackrel{\delta_n}{\rightarrow} \left \{ \begin{array}{cl} {\mathcal B}_{n-1}(F) \otimes F^* &\quad \mbox{if}\, n\geq 3\\ \bigwedge^2F^* & \quad\mbox{if}\, n=2  \end{array} \right . \]

\begin{equation}
\delta_n(\{x\}) = \left\{ \begin{array}{cl}\{x\}_{n-1}\otimes x & \quad \mbox{if}\, n\geq 3\\ (1-x) \wedge x & \quad \mbox{if}\, n=2 \\0 &  \quad \mbox{if}\, \{x\}=\{0\},\{1\},\{\infty\}\end{array}\right .
\end{equation}

Then one defines
\begin{equation}
\mathcal{A}_n(F):=\mathrm{ker}\, \delta_n.
\end{equation}

Note that any element $\alpha(t) = \sum n_i \{f_i(t)\} \in \Z[\PP^1_{F(t)}]$ has a specialization
$\alpha(t_0) = \sum n_i \{f_i(t_0)\} \in \Z[\PP^1_{F}]$, for every $t_0 \in \PP^1_F$.

Thus,
\begin{equation}
\mathcal{R}_n(F) : = \left <  \alpha(0) - \alpha(1); \alpha(t) \in \mathcal{A}_n(F(t))\right >.
\end{equation}

Goncharov proves that $\mathcal{R}_n(\C)$ is the subgroup of all the rational functional equations for the $n$-polylogarithm in $\C$. As stated before, the philosophy is that $\mathcal{R}_n(F)$ should be the subgroup of all the rational functional equations for the $n$-polylogarithm in $F$.

Because of $\delta_n(\mathcal{R}_n(F)) = 0$, it induces morphisms in the quotients
\[ \delta_n :{\mathcal B}_n(F) \rightarrow  {\mathcal B}_{n-1}(F)\otimes F^* \quad n \geq 3, \quad   \delta_2 : {\mathcal B}_2(F) \rightarrow  \bigwedge^2 F^*. \]
One obtains the complex:
\[ {\mathcal B}_F(n): {\mathcal B}_n(F) \stackrel{\delta}{\rightarrow} {\mathcal B}_{n-1}(F)\otimes F^* \stackrel{\delta}{\rightarrow}{\mathcal B}_{n-2}(F)\otimes \bigwedge^2 F^*\stackrel{\delta}{\rightarrow} \dots \stackrel{\delta}{\rightarrow}{\mathcal B}_{2}(F)\otimes \bigwedge^{n-2} F^*\stackrel{\delta}{\rightarrow} \bigwedge^n F^*\]
where
\[ \delta: \{x\}_p\otimes \bigwedge_{i=1}^{n-p}y_i \rightarrow \delta_p(\{x\}_p) \wedge \bigwedge_{i=1}^{n-p}y_i.\]

The following conjecture relates the cohomology of the complex  ${\mathcal B}_F(n)$ to motivic cohomology:
\begin{conj} \label{gonconj} \cite{G4}
\begin{equation}
H^i({\mathcal B}_F(n) \otimes \Q ) \cong  gr^\gamma_n K_{2n-i}(F)\otimes \Q .
\end{equation}
\end{conj}
Evidence supporting this conjecture is found, for instance, in the cases $n=1,2$. First, it is clear that $H^1(\mathcal{B}_F(1))  \cong F^*=  K_1(F)$.

For $n=2$ it is known that
\[\mathcal{B}_2(F) \cong \Z[\PP^1_F]/\left< R_2(x,y); x,y \in F^*, \{0\}, \{\infty\} \right>\]
where
\[ R_2(x,y):=\{x\}+\{y\}+\{1-xy\}+\left\{\frac{1-x}{1-xy}\right\} + \left\{\frac{1-y}{1-xy}\right\} \]
is the five-term relation of the dilogarithm.

Besides,
\begin{equation}\label{sus}
H^1(\mathcal{B}_F(2))_\Q \cong K_3^{\mathrm{ind}}(F)_{\Q}
\end{equation}
\begin{equation}
H^2(\mathcal{B}_F(2)) \cong  K_2(F)
\end{equation}
\begin{equation}
H^n(\mathcal{B}_F(n))  \cong  K_n^M(F)
\end{equation}
The first assertion was proved by Suslin. The second one is Matsumoto's theorem, and the last one corresponds to the definition of Milnor's $K$-theory.

\subsection{Regulators}

Deninger \cite{D} observed that the Mahler measure can be seen as a regulator evaluated in a cycle that may or may not have trivial boundary. More precisely,
\begin{equation}
m(P)  = m(P^*) + \frac{1}{(-2\pi\ii)^{n-1}} \int_G \eta_n(n)(x_1,\dots,x_n).
\end{equation}
We have to explain the ingredients in this formula. In this section we will be concerned with $\eta_n(n)(x_1,\dots,x_n)$. This form will be described in the context of Goncharov's construction of the regulator on the polylogarithmic motivic complexes.

Let us establish some notation:
\[\widehat{\mathcal{L}}_n (z) : = \left \{\begin{array}{cl}
\mathcal{L}_n (z) & n >1 \,\, \mbox{odd} \\ \ii \mathcal{L}_n (z) & n\,\, \mbox{even}\end{array}\right .\]

For any integers $p\geq 1$ and $k\geq 0$, define
\[\beta_{k,p} : = (-1)^p\frac{(p-1)!}{(k+p+1)!} \sum_{j=0}^{\left[\frac{p-1}{2}\right]} \binom{k+p+1}{2j+1} 2^{k+p-2j}B_{k+p-2j}\]
where the $B_i$ are Bernoulli numbers.
\begin{defn}
\[\widehat{\mathcal{L}}_{p,q}(x) : = \widehat{\mathcal{L}}_p(x) \log^{q-1}|x| \dd \log|x| \qquad p\geq 2 \]
\[ \widehat{\mathcal{L}}_{1,q}(x) := (\log|x|\dd \log|1-x|- \log|1-x| \dd \log|x|) \log^{q-1}|x|\]
\end{defn}

Recall that
\[ \mathrm{Alt}_m F(t_1,\dots t_m) : =\sum_{\sigma \in S_m} (-1)^{|\sigma|} F(x_{\sigma(1)}, \dots,  x_{\sigma(m)}).\]

Now, we are ready to describe the differential forms:
\begin{defn} Let $x$, $x_i$ rational functions on a complex variety $X$.
\[\eta_{n+m}(m+1) : \{x\}_n \otimes x_1 \wedge \dots \wedge x_m \rightarrow\]
\[ \widehat{\mathcal{L}}_n(x) \mathrm{Alt}_m \left(\sum_{p\geq0}\frac{1}{(2p+1)!(m-2p)!} \bigwedge_{j=1}^{2p} \dd \log |x_j| \wedge \bigwedge_{j=2p+1}^m  \dd \ii \arg x_j \right) \]
\begin{equation}
+ \sum_{1 \leq k,\, 1\leq p \leq m} \beta_{k,p} \widehat{\mathcal{L}}_{n-k,k}(x) \wedge \mathrm{Alt}_m \left ( \frac{\log|x_1|}{(p-1)!(m-p)!}\bigwedge_{j=2}^{p} \dd \log |x_j| \wedge \bigwedge_{j=p+1}^m \dd \ii\arg x_j   \right)
\end{equation}

\[\eta_{m}(m) : x_1 \wedge \dots \wedge x_m \rightarrow\]
\begin{equation}
\mathrm{Alt}_m \left( \sum_{p\geq0}\frac{\log|x_1|}{(2p+1)!(m-2p-1)!} \bigwedge_{j=2}^{2p+1} \dd \log |x_j| \wedge \bigwedge_{j=2p+2}^m  \dd \ii \arg x_j \right)
\end{equation}

\end{defn}

These differential forms typically will have singularities. In order to work with them we need to have control of the residues. Let $F$ be a field with discrete valuation $v$, residue field $F_v$, and group of units $U$. Let $u \rightarrow \bar{u}$ the projection $U \rightarrow F_v^*$, and $\pi$ a uniformizer for $v$. There is a homomorphism
\[ \theta: \bigwedge^n F^* \rightarrow \bigwedge^{n-1}F_v^*\]
defined by
\[ \theta(\pi\wedge u_1\wedge \dots \wedge u_{n-1})= \bar{u}_1 \wedge \dots \wedge  \bar{u}_{n-1} \qquad \theta(u_1\wedge \dots \wedge u_n) =0.\]

Now define $s_v: \Z[\PP^1_{F}] \rightarrow \Z[\PP^1_{F_v}] $ by $s_v(\{x\}) = \{\bar{x}\}$. It induces $s_v: {\mathcal B}_m(F) \rightarrow {\mathcal B}_m(F_v)$.
Then
\begin{equation}
\partial_v:=s_v\otimes\theta: {\mathcal B}_m(F)\otimes\bigwedge^{n-m}F^* \rightarrow {\mathcal B}_m(F_v)\otimes\bigwedge^{n-m-1}F_v^*
\end{equation}
defines a morphism of complexes
\begin{equation}
\partial_v: {\mathcal B}_F(n) \rightarrow {\mathcal B}_{F_v}(n-1)[-1] .
\end{equation}
\begin{obs} The induced morphism
\[ \partial_v: H^n( {\mathcal B}_F(n)) \rightarrow H^{n-1}( {\mathcal B}_{F_v}(n-1))\]
coincides with the tame symbol defined by Milnor
\[ \partial_v: K^M_n(F) \rightarrow K^M_{n-1}(F_v).\]
\end{obs}

Let $X$ be a complex variety. Let $X^{(1)}$ denote the set of the codimension one closed irreducible subvarieties. Let $\mathcal{A}^j(X)(k)$ denote the space of smooth $j$-forms with values in $(2\pi \ii)^k \R$. Let $\dd$ be the de Rham differential on $\mathcal{A}^j(X)$ and let $\mathcal{D}$ be the de Rham differential on distributions. So
\[ \dd(\dd \arg x) =0 \qquad \mathcal{D}(\dd \arg x) = 2 \pi \delta(x) \]
The difference $\mathcal{D}-\dd$ is the de Rham residue homomorphism.

Goncharov \cite{G4} proves the following,
\begin{thm} \label{thm:cond}  $\eta_n(m)$ induces a homomorphism of complexes
\[ \begin{array}{ccccccc} {\mathcal B}_n(\C(X)) & \stackrel{\delta}{\rightarrow} &{\mathcal B}_{n-1}(\C(X)) \otimes \C(X)^* & \stackrel{\delta}{\rightarrow} & \dots &  \stackrel{\delta}{\rightarrow} & \bigwedge^n\C(X)^*\\\\\downarrow \eta_n(1) & &\downarrow \eta_n(2) &&&&\downarrow \eta_n(n)\\\\\mathcal{A}^0(X)(n-1) & \stackrel{\dd}{\rightarrow} & \mathcal{A}^1(X)(n-1) & \stackrel{\dd}{\rightarrow} & \dots & \stackrel{\dd}{\rightarrow} & \mathcal{A}^{n-1}(X)(n-1) \end{array} \]
such that
\begin{itemize}
\item $\eta_n(1)(\{x\}_n)=\widehat{\mathcal{L}}_n(x)$.
\item $\dd \eta_n(n)(x_1\wedge \dots \wedge x_n) = \beta \re_n\left(\frac{\dd x_1}{x_1} \wedge \dots \wedge \frac{\dd x_n}{x_n}\right)$, where $\beta = 1$ if $n$ is odd and $\ii$ if $n$ is even.
\item $\eta_n(m)(*)$ defines a distribution on $X(\C)$.
\item The morphism $\eta_n(m)$ is compatible with residues:
\begin{equation}
\mathcal{D}\circ \eta_n(m) - \eta_n(m+1) \circ \delta = 2 \pi \ii \sum_{Y \in X^{(1)}} \eta_{n-1}(m-1)\circ \partial_{v_Y}, \qquad m < n
\end{equation}
\begin{equation}
\mathcal{D}\circ \eta_n(n) - \beta {\re}_n \left(\frac{\dd x_1}{x_1} \wedge \dots \wedge \frac{\dd x_n}{x_n}\right) = 2 \pi \ii \sum_{Y \in X^{(1)}} \eta_{n-1}(n-1)\circ \partial_{v_Y},
\end{equation}
where $v_Y$ is the valuation defined by the divisor $Y$.
\end{itemize}
\end{thm}

The relation of $\eta_n(\cdot)$ to the regulator is roughly as follows. As we mentioned in Conjecture \ref{gonconj}, the cohomology of the first complex corresponds to the Adams filtration which is the absolute cohomology. On the other hand, a slight modification of the second complex leads to Deligne cohomology. Now, $\eta_n(\cdot)$ as seen as a map between the cohomologies of these two complexes is conjectured to have the same image as the regulator (\cite{G5}).

We should also remark that the final goal is not to work with $\C(X)$, but with $X$ itself. For $X$ a regular projective variety over a field $F$ Goncharov \cite{G5} describes the difficulties for defining the complex $B_X(n)$ as opposed to $B_{F(X)}(n)$. Basically, it is known how to define $B_X(n)$ when $X=Spec(F)$ for an arbitrary field $F$, $X$ is a regular curve over an arbitrary field $F$, or $X$ is an arbitrary regular scheme but $n\leq 3$. Now, in terms of the relation between $\eta_n(\cdot)$, the regulator, and Beilinson's  conjectures, the picture is much less known. As an illustration, the case $X=Spec(F)$ for $F$ a number field corresponds to Zagier's conjecture \cite{Z1}, \cite{GZ}.

\section{The two-variable case}\label{2-variables}

Rodriguez-Villegas \cite{RV} has performed the explicit construction of the regulator and applied the ideas of Deninger \cite{D} to explain many examples in two variables. This work was later continued by Boyd and Rodriguez-Villegas \cite{BRV1, BRV2}.

Let $P \in \C[x,y]$. Then we may write
\[P(x,y) = a_d(x)y^d + \dots + a_0(x)\]
\[P(x,y) = a_d(x) \prod_{n=1}^d (y-\alpha_n(x)).\]
By Jensen's formula,
\begin{equation}\label{eq:Jensen}
m(P) = m(a_d) + \frac{1}{2 \pi \ii} \sum_{n=1}^d \int_{\TT^1} \log^+ | \alpha_n(x)|
\frac{\dd x}{x} = m(P^*) - \frac{1}{2 \pi } \int_\gamma \eta(x,y).
\end{equation}
Here 
\[ \eta(x,y) := -\ii \eta_2(2)(x \wedge y)= \log|x| \dd  \arg y -\log|y| \dd \arg x\]
is the regulator in this case, defined in the set $C = \{P(x,y) = 0\}$ minus the set $Z$ of zeros ans poles of $x$ and $y$. Also, $P^*=a_d(x)$, and $\gamma$ is the union of paths in $C$ where
$|x| = 1 $ and $|y| \geq 1$. Finally, note that $\partial \gamma = \{ P(x,y)=0 \}\cap \{ |x|=|y|=1\}$.

In general, $\eta(x,y)$ is closed in $C \setminus Z$, since $\dd \eta(x,y) = \im \left( \frac{\dd x}{x}\wedge \frac{\dd y}{y}\right)$ (see Theorem \ref{thm:cond}). Now, if we would like to be able to perform this computation, we wish to arrive to one of these two situations:
\begin{enumerate}
\item $\eta$ is exact, and $\partial \gamma \not = 0$. In this case we
can integrate using Stokes Theorem.
\item $\eta$ is not exact and $\partial \gamma = 0$. In this case we
can compute the integral by using the Residue Theorem.
\end{enumerate}

Examples for the second case are found, for instance, in the family of Laurent polynomials $x+\frac{1}{x}+y+\frac{1}{y}+k$ studied by Boyd \cite{B2}, Deninger \cite{D}, and Rodriguez-Villegas \cite{RV}. Technically one needs $k \not \in [-4,4]$ for these examples to be in the first case, otherwise  $\partial \gamma \not = 0$. However Deninger has given an interpretation that allows an adaptation of the cases of $k\in (-4,4) \setminus \{0\}$ into this frame as well.

Following Theorem \ref{thm:cond},
\begin{equation}
\eta(x,1-x)=\dd D(x).
\end{equation}
Thus, $\eta$ is exact when
\begin{equation}\label{cond0}
x\wedge y = \sum_i r_i x_i\wedge (1-x_i)
\end{equation}
in $\bigwedge^2(\C(C)^*)\otimes \Q$. This condition may be rephrased as the symbol $\{x,y\}$ is trivial in $K_2(\C(C))$.

In fact, if condition (\ref {cond0}) is satisfied, we obtain
\begin{equation}
\eta(x,y) = \sum_i r_i \dd D(x_i)  = \dd D\left(\sum_i r_i \{x_i\}_2\right).
\end{equation}
Finally we write
\[ \partial \gamma = \sum_k \epsilon_k [w_k], \quad \epsilon_k = \pm 1\]
where $w_k \in C(\C)$, $|x(w_k)| = |y(w_k)| =1$. Thus
\begin{thm}\cite{BRV1,BRV2}
Let $P\in \C[x,y]$ be irreducible and such that $x,y$ satisfies equation (\ref{cond0}). Then
\[2 \pi (\log|a_d|-m(P)) = D(\xi) \qquad \mbox{for} \quad \xi =\sum_k \sum_i \epsilon_k r_i \{x_i(w_k)\}_2. \]
\end{thm}
Boyd and Rodriguez-Villegas prove even more. Under certain assumptions, it is possible to apply Zagier's Theorem \cite{Z0} and relate the Mahler measure of $P$ to a rational combination of terms of the form $\frac{|\Delta|^{\frac{1}{2}}\zeta_F(2)}{\pi^{2[F:\Q]-2}}$ for certain number fields $F$ which depend on $P$ (or more specifically, on $w_k$).

\subsection{An example for the two-variable case}
\begin{figure}
\centering
\includegraphics[width=12pc]{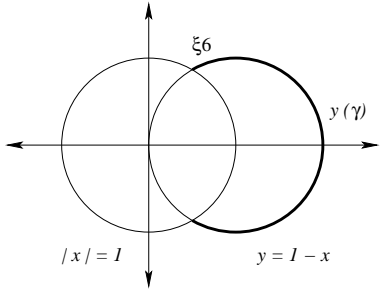}
\caption{\label{figsmyth} Integration path for $x+y-1$}
\end{figure}

To be concrete, we are going to examine the simplest example for the exact case in two variables. Consider Smyth's formula:
\[\pi m(x+y-1) = \frac{3 \sqrt{3}}{4} \Lf (\chi_{-3}, 2).\]
For this case,
\[ x \wedge y = x \wedge (1-x).\]
Then
\[ 2 \pi m(P) = - \int_\gamma \eta(x,y)= - \int_\gamma \eta(x,1-x) = -D(\partial \gamma)\]
Here
\[\gamma=\{(x,y) \, | \, |x|=1, |1-x|\geq 1\} = \{(\e^{2 \pi \ii \theta}, 1-\e^{2 \pi \ii \theta})\, | \,    \theta \in \left[ 1 / 6\, ;  5 / 6 \right] \}.\]
Figure \ref{figsmyth} shows the integration path $\gamma$.

Then  $\partial \gamma = [\bar{\xi}_6] - [{\xi}_6]$ (where $\xi_6 = \frac{1+\sqrt{3} \ii}{2}$) and we obtain
\[2 \pi m(x+y-1) = D(\xi_6) - D(\bar{\xi_6}) = 2  D(\xi_6) = \frac{3 \sqrt{3}}{2} \Lf (\chi_{-3}, 2).\]

\section{The three-variable case}\label{3-var}

Our goal is to extend this situation to three variables. Let $P \in \C[x,y,z]$. We will take

\[\eta(x,y,z) := \eta_3(3)(x\wedge y \wedge z)= \log|x|\left( \frac{1}{3} \dd \log |y| \wedge \dd \log |z| - \dd \arg y \wedge \dd \arg z\right )\]
\[ + \log|y|\left( \frac{1}{3} \dd \log |z| \wedge \dd \log |x| - \dd \arg z \wedge \dd \arg x\right)\]
\begin{equation}
 +  \log|z|\left( \frac{1}{3} \dd \log |x| \wedge \dd \log |y| - \dd \arg x \wedge \dd \arg y\right).
\end{equation}
This differential form is defined in the surface $S = \{ P(x,y,z) = 0 \}$ minus the set $Z$ of poles and zeros of $x$, $y$ and $z$.

We can express the Mahler measure of $P$ as
\begin{equation} \label{3var}
m(P) = m(P^*)- \frac{1}{(2 \pi )^2} \int_{\Gamma} \eta(x,y,z).
\end{equation}
Where $P^*$, following the previous notation, is the principal coefficient of the polynomial $P \in \C[x,y][z]$ and
\[\Gamma = \{ P(x,y,z) = 0 \} \cap \{ |x|=|y| =1, |z| \geq 1 \}.\]

Recall $\eta$ in closed in $S \setminus Z$ since it verifies $\dd \eta(x,y,z) = \re \left( \frac{\dd x}{x} \wedge \frac{\dd  y}{y} \wedge \frac{\dd z}{z}\right)$.
Typically, one expects that integral (\ref{3var}) can be computed if we are in one of the two ideal situations that we described before. Either the form $\eta(x,y,z)$ is not exact and the set $\Gamma$ consists of closed subsets and the integral is computed by residues, or the form $\eta(x,y,z)$ is exact and the set $\Gamma$ has nontrivial boundaries, so Stokes Theorem is used.

The first case would lead to instances of Beilinson's conjectures and produces special values of L-functions of surfaces. Examples in this direction can be found in Bertin's work \cite{Be}. Bertin relates the Mahler measure of some $K3$ surfaces to Eisenstein-Kronecker series in a similar way as Rodriguez-Villegas does for two-variable cases \cite{RV}.

In the second case we need that $\eta(x,y,z)$ is exact. We are going to concentrate on this case.

We are integrating on a subset of the surface $S$. In order for the element in the cohomology to be defined everywhere in the surface $S$, we need the residues to be zero. This situation is fulfilled when the tame symbols are zero (see Section \ref{Gonchcomplexes}). This condition will not be a problem for us because  when $\eta$ is exact the tame symbols are zero.

As in the two-variable case, Theorem \ref{thm:cond} implies
\begin{equation}
\eta (x, 1-x, y) = \dd \,\omega (x,y)
\end{equation}
where
\begin{equation}
 \omega(x,y)  := \eta_3(2)(\{x\}_2\otimes y) = - D(x) \dd \arg y + \frac{1}{3}\log |y| ( \log |1-x| \dd \log |x| - \log |x| \dd \log |1-x|).
\end{equation}

Thus, in order to apply Stokes Theorem, we need to require that
\begin{equation}
x \wedge y \wedge z = \sum r_i x_i \wedge (1-x_i) \wedge y_i
\end{equation}
in $\bigwedge^3(\C(S)^*)\otimes \Q$ for $\eta$ to be exact. An equivalent way of expressing this condition is that
$\{x,y,z\}$ is trivial in $K^M_3(\C(S))$.

In this case,
\[\int_{\Gamma} \eta(x,y,z)  = \sum r_i  \int_{\Gamma} \eta(x_i , 1-x_i ,y_i) = \sum r_i  \int_{\partial \Gamma} \omega (x_i , y_i),    \]
where
\[\partial \Gamma = \{ P(x,y,z) = 0 \} \cap \{ |x|=|y| =|z| = 1 \}.\]
This set $\partial \Gamma$ seems to have no boundary. However, $\partial \Gamma$ as described above may contain singularities which may give rise to a boundary when desingularized. We will change our point of view. Namely, assume that $P\in \R[x,y,z]$ and nonreciprocal (this condition is true for all the examples we study), then
\[P(x,y,z)=P(\bar{x},\bar{y},\bar{z}).\]
This property, together with the condition $|x|=|y|=|z|=1$, allows us to write
\[\partial \Gamma = \{ P(x,y,z) =  P(x^{-1},y^{-1}, z^{-1}) = 0  \} \cap \{ |x|=|y| =1 \}. \]
(This idea was proposed by Maillot). Observe that we are integrating now on a path $\{ |x|=|y| =1 \}$ inside the curve
\[ C = \{ \res_z(P(x,y,z), P(x^{-1},y^{-1}, z^{-1})) = 0  \}.\]

In order to easily compute
\[ \int_{\partial \Gamma} \omega(x,y)\]
we have again the two possibilities that we had before. We are going to concentrate, as usual, in the case when $\omega(x,y)$ is exact.

The differential form $\omega$ is defined in this new curve $C$. As before, to be sure that it is defined everywhere, we need to ask that the residues are trivial. This fact is guaranteed by the triviality of tame symbols. This last condition is satisfied if $\omega$ is exact. Indeed, we have changed our ambient variety, and we now wonder when $\omega$ is exact in $C$ ($\omega$ is not exact in $S$ since that would imply that $\eta$ is zero).


Fortunately we have
\begin{equation}
 \omega(x,x) = \dd \mathcal{L}_3(x)
\end{equation}
by Theorem \ref{thm:cond}.

The condition for $\omega$ to be exact  is not as easily established as in the preceding cases because $\omega$ is not multiplicative in the first variable. In fact, the first variable behaves as the dilogarithm, in other words, the transformations are ruled by the five-term relation.
We may express the condition we need as:
\begin{equation} 
\{x\}_2 \otimes y  = \sum r_i \{x_i\}_2 \otimes x_i
\end{equation}
in $(\mathcal{B}_2(\C(C)) \otimes \C(C)^*)_\Q$. Assuming Conjecture \ref{gonconj}, this is equivalent to saying that a certain symbol for $x$ and $y$ is trivial in  $gr^\gamma_3 K_{4}(\C(C))\otimes \Q$ .
Then we have
\[\int_{\gamma} \omega (x , y)  =  \sum r_i  \left. \mathcal{L}_3(x_i) \right|_{\partial \gamma}.\]
where $\gamma = C \cap \TT^2$.

Now assume that
\[ \partial \gamma = \sum_k \epsilon_k [w_k], \quad \epsilon_k = \pm 1\]
where $w_k \in C(\C)$, $|x(w_k)| = |y(w_k)| =1$. Thus we have proved

\begin{thm} Let $P(x,y,z) \in \R[x,y,z]$ be irreducible and nonreciprocal, and let $S=\{P(x,y,z)=0\}$ and $C=\{ \res_z(P(x,y,z), P(x^{-1},y^{-1}, z^{-1})) = 0 \}$.
Assume that
\begin{equation}\label{eqcond1}
x \wedge y \wedge z = \sum_i r_i x_i \wedge (1-x_i) \wedge y_i
\end{equation}
in $\bigwedge^3(\C(S)^*)\otimes \Q$, and
\begin{equation}\label{eqcond2}
\{x_i\}_2 \otimes y_i  = \sum_{j} r_{i,j} \{x_{i,j}\}_2 \otimes x_{i,j}
\end{equation}
in $(\mathcal{B}_2(\C(C)) \otimes \C(C)^*)_\Q$ for all $i$. Then
\begin{equation}
4\pi^2 (m(P^*)-m(P)) = \mathcal{L}_3(\xi)\qquad \mbox{for} \quad \xi =\sum_k \sum_{i, j} \epsilon_k r_i r_{i,j} \{x_{i, j}(w_k)\}_3.
 \end{equation}
\end{thm}

By using Zagier's conjecture (see Zagier \cite{Z1}, Zagier and Gangl \cite{GZ}), it is possible to formulate a conjecture that would imply, under certain additional circumstances, a relationship with $\zeta_F(3)$ in a similar fashion as Boyd and Rodriguez-Villegas have done  for the two-variable case. We will illustrate this phenomenon at the end of the example that follows.

\subsection{The case of $\res_{\{0,m,m+n\}}$}

We will proceed to the study of a family of three-variable polynomials that come from the world of resultants, namely, $\res_{\{0,m,m+n\}}$. This family  was computed in \cite{DL} and the computation is quite involved, though elementary. The Mahler measure of $\res_{\{0,m,m+n\}}$ is the same as the Mahler measure of a certain rational function. More precisely,
\begin{thm}\cite{DL}
\begin{equation}
m\left(z - \frac{(1-x)^m (1-y)^n}{(1-xy)^{m+n}}\right)=\frac{2n}{\pi^2} (\mathcal{L}_3(\phi_2^{m}) - \mathcal{L}_3(-\phi_1^{m})) +\frac{2m}{\pi^2} (\mathcal{L}_3(\phi_1^{n}) - \mathcal{L}_3(-\phi_2^{n}))
\end{equation}
where $\phi_1$ is the root of $x^{m+n}+x^n-1=0$ that lies in the interval $[0,1]$ and $\phi_2$ is the root of $x^{m+n}-x^n-1=0$ that lies in $[1, \infty)$.
\end{thm}
\begin{proof}
Since we would like to see that $\eta(x,y,z)$ is exact, we need to solve equation (\ref{eqcond1}) for this case. The equation for the wedge product becomes
\begin{eqnarray*}
x \wedge y \wedge z & = & m x \wedge y \wedge (1-x) + n x \wedge y \wedge (1-y)- (m+n) x \wedge y \wedge (1-xy) \\
& = & -m x \wedge (1-x) \wedge y  + n y \wedge (1-y) \wedge x \\
&&+ m xy \wedge (1-xy) \wedge y - n xy \wedge (1-xy) \wedge x.
\end{eqnarray*}

After performing Stokes Theorem for the first time we will have to evaluate the form $\omega$ in the following element of $\mathcal{B}_2(\C(C))\otimes\C(C)^*$:
\[ \Delta = m(  \{xy\}_2 \otimes y - \{x\}_2\otimes y ) -n ( \{xy\}_2 \otimes x  -  \{y\}_2 \otimes x).\]

We need to compute the corresponding curve $C$. We take advantage of the fact that our equation has the shape $z = R(x,y)$. In order to compute $C$, we simply need to consider
\begin{equation}\label{eq:trick}
R(x,y)R(x^{-1},y^{-1}) = z \cdot z^{-1} = 1.
\end{equation}
For this case
\[ \frac{(1-x)^m (1-y)^n (1-x^{-1})^m (1-y^{-1})^n}{(1-xy)^{m+n}(1-x^{-1}y^{-1})^{m+n}} = 1. \]

Let us denote
\[ x_1 = \frac{1-x}{1-xy} \quad y_1 = \frac{1-y}{1-xy} \quad \widehat{x}_1 = 1 - x_1 \quad \widehat{y}_1 = 1 - y_1. \]

then we may rewrite the equation for $C$ as
\[x_1^m y_1^n  \widehat{x}_1^n  \widehat{y}_1^m = 1. \]

Now we use the five-term relation:
\[ \{x\}_2 + \{y\}_2 + \{1-xy\}_2 + \{ x_1\}_2 + \{ y_1\}_2 = 0.\]

Then we obtain
\[ \Delta = m ( \{y\}_2 \otimes y+ \{x_1\}_2 \otimes y + \{y_1\}_2 \otimes y  )- n ( \{x\}_2 \otimes x+ \{x_1\}_2 \otimes x + \{y_1\}_2 \otimes x   ). \]

Observe that $x = \frac{\widehat{x}_1}{y_1} $, $y = \frac{\widehat{y}_1}{x_1} $.

Thus, we may write
\[\Delta= m ( \{y\}_2 \otimes y+ \{x_1\}_2 \otimes \widehat{y}_1 - \{x_1\}_2 \otimes x_1  +  \{y_1\}_2 \otimes \widehat{y}_1 - \{y_1\}_2 \otimes x_1    )   \]
\[ - n ( \{x\}_2 \otimes x+ \{x_1\}_2 \otimes \widehat{x}_1 - \{x_1\}_2 \otimes y_1  +  \{y_1\}_2 \otimes \widehat{x}_1 - \{y_1\}_2 \otimes y_1    )   \]
\[= m \{y\}_2 \otimes y+  \{x_1\}_2 \otimes \widehat{y}^m_1 - m \{x_1\}_2 \otimes x_1  -  m \{\widehat{y}_1\}_2 \otimes \widehat{y}_1 -  \{y_1\}_2 \otimes x^m_1      \]
\[ - n  \{x\}_2 \otimes x+ n \{\widehat{x}_1\}_2 \otimes \widehat{x}_1 + \{x_1\}_2 \otimes y^n_1  -  \{y_1\}_2 \otimes \widehat{x}^n_1 + n \{y_1\}_2 \otimes y_1.       \]

Because of the equation for $C$,
\[ \{x_1\}_2 \otimes y_1^n \widehat{y}_1^m   - \{y_1\}_2 \otimes x_1^m \widehat{x}_1^n  = - \{x_1\}_2 \otimes x_1^m \widehat{x}_1^n + \{y_1\}_2 \otimes y_1^n \widehat{y}_1^m \]
\[= - m \{x_1\}_2 \otimes x_1 + n  \{\widehat{x}_1\}_2 \otimes \widehat{x}_1 + n  \{y_1\}_2 \otimes y_1 - m  \{\widehat{y}_1\}_2 \otimes \widehat{y}_1,  \]
we obtain,
\[ \Delta = m ( \{y\}_2 \otimes y- \{\widehat{y}_1\}_2 \otimes \widehat{y}_1 - \{x_1\}_2 \otimes x_1  -  \{\widehat{y}_1\}_2 \otimes \widehat{y}_1 - \{x_1\}_2 \otimes x_1    )   \]
\[ - n ( \{x\}_2 \otimes x- \{\widehat{x}_1\}_2 \otimes \widehat{x}_1 - \{y_1\}_2 \otimes y_1  -  \{\widehat{x}_1\}_2 \otimes \widehat{x}_1 - \{y_1\}_2 \otimes y_1    ) .\]

\[ \Delta = m ( \{y\}_2 \otimes y - 2 \{\widehat{y}_1\}_2 \otimes \widehat{y}_1 - 2 \{x_1\}_2 \otimes x_1 )   - n ( \{x\}_2 \otimes x - 2\{\widehat{x}_1\}_2 \otimes \widehat{x}_1 - 2 \{y_1\}_2 \otimes y_1 )   \]

We now need to study the path of integration. First write $x = \e^{ 2\ii \alpha}$, $ y = \e^{2\ii \beta}$, for $ -\frac{\pi}{2}  \leq \alpha, \beta \leq \frac{\pi}{2}$. Then,
\[ x_1 = \e^{-\ii \beta} \frac{\sin \alpha}{\sin (\alpha+\beta)} \quad y_1 = \e^{-\ii \alpha} \frac{\sin \beta}{\sin (\alpha+\beta)} \]
and\[ \widehat{x}_1 = \e^{\ii \alpha} \frac{\sin \beta}{\sin (\alpha+\beta)} \quad \widehat{y}_1 = \e^{\ii \beta} \frac{\sin \alpha}{\sin (\alpha+\beta)}.\]

Let $a = \left|\frac{\sin \alpha}{\sin (\alpha + \beta)}\right|$, $b =  \left| \frac{\sin \beta}{\sin (\alpha + \beta)}\right|$. Then we may write
\[ x_1 = \pm a \e^{-\ii \beta} \quad y_1 = \pm b \e^{-\ii \alpha} \quad \widehat{x}_1 =\pm b \e^{\ii \alpha} \quad \widehat{y}_1 =\pm a \e^{\ii \beta}. \]

By means of the Sine theorem, we may think of $a$, $b$ and 1 as the sides of a triangle with the additional condition
\[ a^m b^n = 1.\]
The triangle determines the angles, $\alpha$ and $\beta$, which are opposite to the sides $a$, $b$ respectively. We need to be careful and take the complement of an angle if it happens to be greater than $\frac{\pi}{2}$, (this corresponds to the cases when the sines are negatives). However, we need to be cautious. In fact, the problem of constructing the triangle given the sides has always two symmetric solutions. We are going to count each triangle once, so we will need to multiply our final result by two. To sum up, $a$ and $b$ are enough to describe the set where the integration is performed.

\begin{figure}
\centering
\includegraphics[width=28pc]{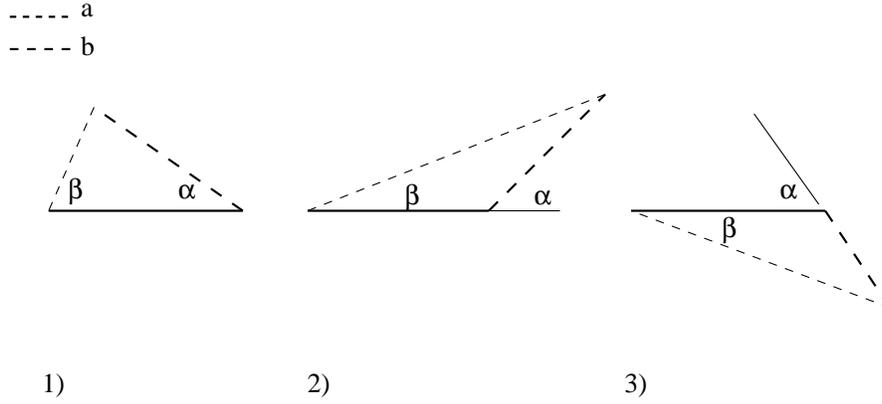}
\caption{\label{mat3res} We are integrating over all the possible triangles. The angles have to be measured negatively if they are greater than $\frac{\pi}{2}$ as $\alpha$ in the case 2). We will not count the triangles pointing down as in 3). }
\end{figure}

Now, the boundaries (where the triangle degenerates) are three: $b+1=a$, $a+1=b$ and $a+b=1$. Let
\[ \begin{array}{ccccl}
 \phi_1  & \mbox{be the root of} & x^{m+n}+x^n-1 =0,  &  \mbox{ with} & 0 \leq \phi_1 \leq 1,\\
 \phi_2  & \mbox{be the root of} &   x^{m+n}-x^n-1 =0,  & \mbox{ with} & 1 \leq \phi_2.
\end{array}\]

Then the first two conditions are translated as
\[ \begin{array}{cccc}
  a= \phi_1^{-n},& b= \phi_1^m,  & \alpha=0,& \beta=0,\\
  a= \phi_2^{-n}, &b= \phi_2^m,  & \alpha=0,& \beta=0.
\end{array}\]

The third condition is inconsequential, since it requires both $a, b \leq 1$ (but they can not be both equal to 1 at the same time) and $a^m b^n=1$.

Hence, the integration path (from condition $a+1=b$ to $b+1=a$) is
\[ \begin{array}{cc}
 0 \leq \alpha \leq \theta_1,& 0 \geq \beta \geq -\frac{\pi}{2}, \\\\
 \theta_1 \leq \alpha \leq \frac{\pi}{2}, & \frac{\pi}{2} \geq \beta \geq \theta_2,\\\\
  -\frac{\pi}{2} \leq \alpha \leq 0, & \theta_2 \geq \beta \geq 0. \\\\
\end{array}\]

Here $\theta_1$ is the angle that is opposite to the side $a$ when the triangle is right-angled with hypotenuse $b$ and $\theta_2$ is opposite to $b$ when $a$ is the hypotenuse. We do not need to compute those angles. In fact, we may describe the integration path as either
\[ 0 \leq \alpha \leq \frac{\pi}{2}, \qquad - \frac{\pi}{2} \leq \alpha \leq 0, \]
or
\[0 \geq \beta \geq -\frac{\pi}{2}, \qquad  \frac{\pi}{2} \geq \beta \geq 0. \]
It is appropriate to think of it in this way, because $\{x_1\}_3+\{\widehat{y}_1\}_3$ and $\{\widehat{x}_1\}_3+\{y_1\}_3$ change continuously around the right-angled triangles. Moreover, because of this property, everything reduces to evaluating $\mathcal{L}_3$ in
\[\Omega =m(\{y\}_3-2\{\widehat{y}_1\}_3-2\{x_1\}_3)-n(\{x\}_3-2\{\widehat{x}_1\}_3-2\{y_1\}_3) \] in the cases of $b+1=a$ and $a+1=b$ and computing the difference.

One could have problems when $z$ is zero or has a pole. $z$ is zero for $x=1$ and $y=1$, but these conditions correspond to $\Delta =m\{y\}_2\otimes y $ and  $\Delta =-n\{x\}_2\otimes x $. They  lead to $\Omega=m\{y\}_3$ and $\Omega=-n\{x\}_3$ and integrate to zero when the variables move in the unit circle.

The poles are at $xy=1$, which corresponds to $\Delta = (m-n) \{x\}_2 \otimes x$.
Integrating, we obtain $\Omega = (m-n) \{x\}_3$
which leads to zero when $x$ moves in the unit circle.

We obtain
\[ 4 \pi^2 m(P)=2\left( 4n(\mathcal{L}_3(\phi_2^m) - \mathcal{L}_3(-\phi_1^m) ) + 4m (\mathcal{L}_3(\phi_1^{n}) - \mathcal{L}_3(-\phi_2^{n}))  \right). \]

Finally,
\[ m(P) =  \frac{2n}{\pi^2} (\mathcal{L}_3(\phi_2^{m}) - \mathcal{L}_3(-\phi_1^{m})) +\frac{2m}{\pi^2} (\mathcal{L}_3(\phi_1^{n}) - \mathcal{L}_3(-\phi_2^{n})),\]
so we recover the result of \cite{DL}.
\end{proof}

The case with $m=n=1$ is especially elegant. Here the rational function has the form
\[ z = \frac{(1-x)(1-y)}{(1-xy)^2},\]
and
\[m(P)= \frac{4}{\pi^2}(\mathcal{L}_3(\phi) - \mathcal{L}_3(-\phi))\]
where $\phi^2+\phi-1=0$ and $0\leq \phi \leq1$ (in other words, $\phi= \frac{-1+\sqrt{5}}{2}$).

Moreover, we may use Zagier's conjecture to describe this result in terms of the zeta function of $\Q(\sqrt{5})$. According to the conjecture, $H^1(\mathcal{B}_{\Q(\sqrt{5})}(3))$ has rank 2. We may take $\{\{1\}_3, \{\phi\}_3\}$ as basis. In order to see this we need to check that $\{\phi\}_2\otimes \phi$ is trivial. That is the case because
\[\{\phi\}_2=\{1-\phi^2\}_2=-\{\phi^2\}_2=2\{-\phi\}_2 - 2\{\phi\}_2 = -2\{1+\phi\}_2 - 2\{\phi\}_2 = -2\{\phi^{-1}\}_2 - 2\{\phi\}_2 =0,\]
which implies $\{\phi\}_2\otimes \phi=0$. Then the conjecture predicts
\[ \zeta_{\Q(\sqrt{5})}(3) \sim_{\Q^*} \sqrt{5} \left|\begin{array}{cc}\mathcal{L}_3(\phi) & \mathcal{L}_3(1)\\  \mathcal{L}_3(-\phi^{-1}) & \mathcal{L}_3(1)\end{array}\right| =\sqrt{5} \zeta(3) (\mathcal{L}_3(\phi)-\mathcal{L}_3(-\phi)). \]

Indeed,
\[\zeta_{\Q(\sqrt{5})}(3) = \frac{\zeta(3)}{\sqrt{5}} ( \mathcal{L}_3(\phi) - \mathcal{L}_3(-\phi)). \]

Which allows us to write
\[m(  \res_{\{0,1,2\}}) = \frac{4 \sqrt{5} \zeta_{\Q(\sqrt{5})}(3)}{\pi^2 \zeta(3)}.\]

\section{A few words about the four-variable case}

Unfortunately, we do not have a general systematic method to algebraically describe the successive integration domains in more than three variables. Hence, we can not formulate a precise general result. However, this does not prevent us from using a similar technique for some four-variable cases. In this section we recall the list of differentials in four variables.

The sequence of differentials should be as follows:

\[ \eta(x,y,w,z) := -\ii \eta_4(4)(x,y,w,z)= \frac{1}{4} \left( -\log |z| \im\left ( \frac{\dd x}{x} \wedge \frac{\dd y}{y} \wedge \frac{\dd w}{w} \right ) + \log |w|\im\left ( \frac{\dd x}{x} \wedge \frac{\dd y}{y} \wedge \frac{\dd z}{z} \right )  \right.\]
\[\left. -\log |y| \im\left ( \frac{\dd x}{x} \wedge \frac{\dd w}{w} \wedge \frac{\dd z}{z} \right ) + \log|x| \im\left ( \frac{\dd y}{y} \wedge \frac{\dd w}{w} \wedge \frac{\dd z}{z} \right ) \right. \]
\begin{equation}
\left.  + \eta(x,y,w)\wedge \dd \arg z -\eta(x,y,z)\wedge \dd \arg w+ \eta(x,w,z)\wedge \dd \arg y -\eta(y,w,z)\wedge \dd \arg x \right)
\end{equation}
where $\eta(x,y,z)$ denotes the differential previously defined for three variables.

We have,
\begin{equation}
\eta(x,1-x,y,w) = \dd \omega(x,y,w)
\end{equation}
where
\[\omega(x,y,w) := -\ii \eta_4(3)(x,y,w)= D(x) \left( \frac{1}{3} \dd \log|y| \wedge \dd \log|w| - \dd \arg y \wedge \dd \arg w \right)\]
\begin{equation}
 + \frac{1}{3} \eta(y,w) \wedge \left( \log|x| \dd \log|1-x| - \log|1-x| \dd \log|x| \right).
\end{equation}

Next,
\begin{equation}
\omega(x,x,y) = \dd \mu(x,y)
\end{equation}
with
\begin{equation}
\mu(x,y) := -\ii \eta_4(2)(x,y)= \mathcal{L}_3(x) \dd \arg y - \frac{1}{3} D(x) \log|y| \dd \log|x|.
\end{equation}

Finally,
\begin{equation}
\mu(x,x) = \dd \mathcal{L}_4(x).
\end{equation}

\subsection{An example in four variables}
In spite of the fact that we do not know how to treat the integration domains, we may still be able to do the algebraic integration for some examples of four-variable polynomials. Here is an example.

We will study the case of $\res_{\{(0,0),(1,0),(0,1)\}}$, whose Mahler measure was first computed in \cite{DL}. This is the case of the nine-variable polynomial that is the general $3\times3$ determinant. Because of homogeneities, this Mahler measure problem may be reduced to computing the Mahler measure of a four-variable polynomial. The result is
\begin{thm}\cite{DL}
\begin{equation}
m((1-x)(1-y) - (1-w)(1-z))= \frac{9}{2\pi^2} \zeta(3).
\end{equation}
\end{thm}
\begin{proof}
First we have to solve the equation with the wedge product:
\[ x\wedge y \wedge w \wedge z = - \frac{1}{x} \wedge y  \wedge w \wedge z =
- \frac{1}{x} \wedge y\left(1-\frac{1}{x}\right)  \wedge w \wedge z +  \frac{1}{x} \wedge \left(1-\frac{1}{x} \right)  \wedge w \wedge z.\]

Now the first term on the right-hand side is
\[- \frac{1}{x} \wedge y\left(1-\frac{1}{x}\right)  \wedge w \wedge z =
\frac{x}{w} \wedge \left(y-\frac{y}{x}\right) \wedge w \wedge z \]
\[=\frac{x}{w}\left( 1 - y+\frac{y}{x}\right) \wedge \left(y-\frac{y}{x}\right) \wedge w \wedge z - \left( 1 - y+\frac{y}{x}\right) \wedge \left(y-\frac{y}{x}\right)  \wedge w \wedge z. \]

Next, we use the formula for $z$ as a function of the other variables:
\[\frac{x}{w}\left( 1 - y+\frac{y}{x}\right) \wedge \left(y-\frac{y}{x}\right)  \wedge w \wedge z = \frac{x+y-xy}{w}\wedge \left(y-\frac{y}{x}\right)  \wedge w \wedge \frac{-w+x+y-xy}{w( 1-w)}\]
\[= \frac{x+y-xy}{w}\wedge \left(y-\frac{y}{x}\right) \wedge w \wedge \left(1-\frac{x+y-xy}{w}\right) - \frac{x+y-xy}{w}\wedge \left(y-\frac{y}{x}\right)  \wedge w \wedge ( 1-w).\]

Note that
\[- (x+y-xy)\wedge \left(y-\frac{y}{x}\right)  \wedge w \wedge ( 1-w)\]
\[ = -  \left(1-y+\frac{y}{x}\right)\wedge \left(y-\frac{y}{x}\right)  \wedge w \wedge ( 1-w)- x \wedge\left(y-\frac{y}{x}\right)  \wedge w \wedge ( 1-w).\]

Hence
\[ x\wedge y \wedge w \wedge z = \frac{1}{x} \wedge \left(1 - \frac{1}{x} \right)\wedge w \wedge z +  \left(y-\frac{y}{x}\right) \wedge\left( 1 - y+\frac{y}{x}\right)   \wedge w \wedge z(1-w) \]
\[+ \frac{x+y-xy}{w}\wedge \left(1-\frac{x+y-xy}{w}\right)\wedge \left(y-\frac{y}{x}\right) \wedge w -w \wedge ( 1-w)\wedge x\wedge \left(y-\frac{y}{x}\right).\]

The form $\omega$ will be evaluated in the following element:
\[ \Delta = \left\{ \frac{1}{x} \right\}_2 \otimes w \wedge z +  \left\{y-\frac{y}{x}\right\}_2 \otimes  w \wedge z(1-w)\]
\[ + \left\{\frac{x+y-xy}{w}\right\}_2 \otimes \left(y-\frac{y}{x}\right) \wedge w
- \{w\}_2\otimes x \wedge \left(y-\frac{y}{x}\right)   \]

\[ = -\{x\}_2 \otimes w \wedge z +  \left\{y-\frac{y}{x}\right\}_2 \otimes  w \wedge z(1-w)\]
\[ -\left\{z-\frac{z}{w}\right\}_2 \otimes \left(y-\frac{y}{x}\right) \wedge w
- \{w\}_2\otimes x\wedge \left(y-\frac{y}{x}\right). \]

For applying Stokes Theorem, we still apply the technique that is analogous to the computation of the equation for $C$ in the three-variable case. We can apply the analogue of the equation (\ref{eq:trick}):
\[\left( 1 - \frac{(1-x)(1-y)}{1-w}\right) \left( 1 - \frac{(1-x^{-1})(1-y^{-1})}{1-w^{-1}}\right)=1,\]
which can be simplified as
\[ x=1, \quad y=1, \quad w =x, \quad \mbox{or}\quad w=y.\]

The above conditions correspond to two pyramids in the torus $\TT^3$, as seen in picture \ref{cube}. We will make the computation over the lower pyramid and then multiply the result by 2.

\begin{figure}
\centering
\includegraphics[width=39pc]{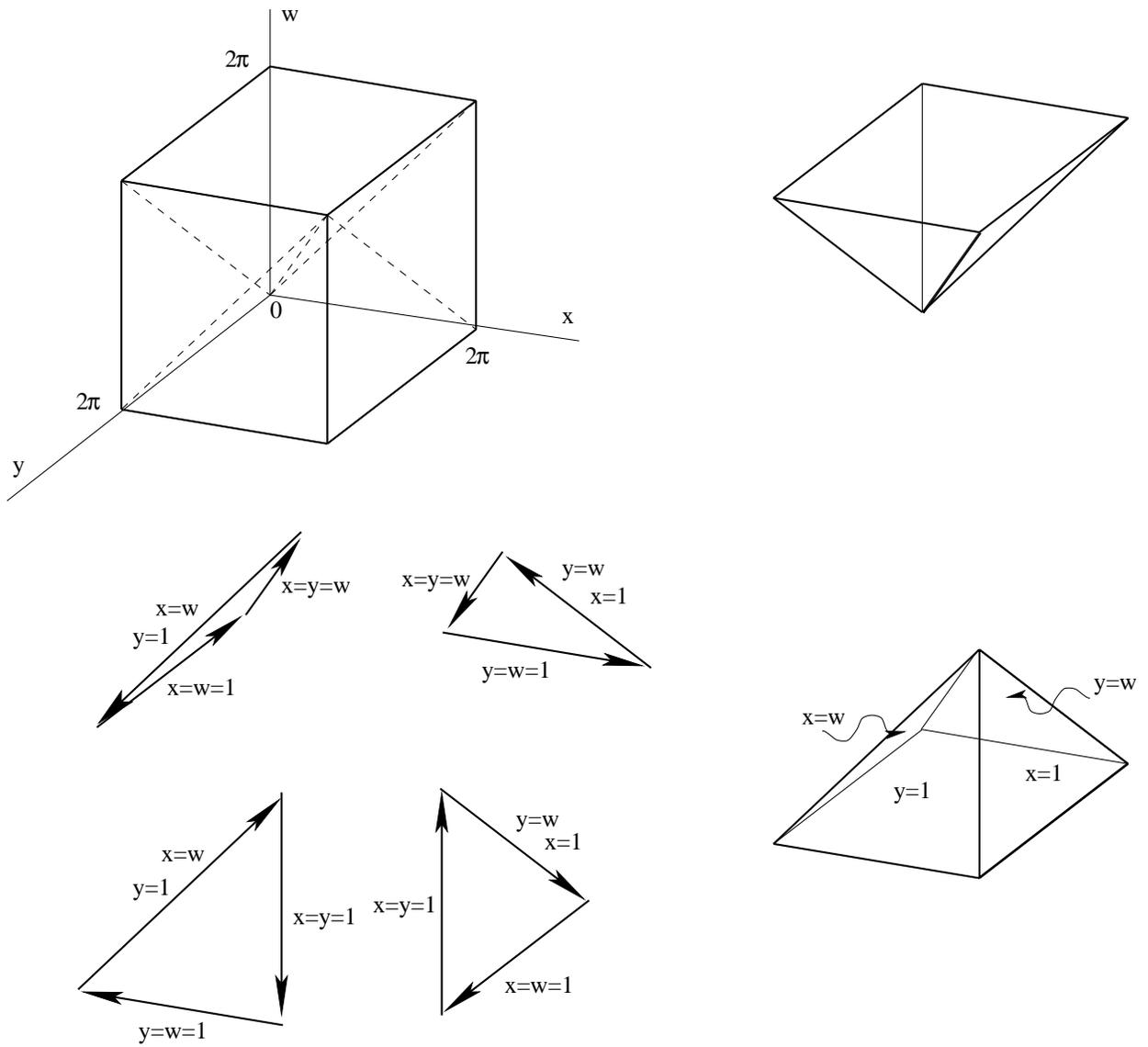}
\caption{\label{cube} Integration set for $\res_{\{(0,0),(1,0),(0,1)\}}$.}
\end{figure}

When $x=1$, in this case, $w=1$ or $z=1$. If $w=1$, $\Delta=0$.

If $z=1$,
\[ \Delta = -\left\{1-\frac{1}{w}\right\}_2\otimes y \wedge w.\]
Then $\mu$ will be evaluated on
\[\Omega= \{w\}_3\otimes y.\]

$\Omega$ will be integrated on the boundary, which is $y=1$, $w=1$ and $y=w$.

If $y=1$, $\Omega=0$.
If $w=1$,
\[\Omega = \{1\}_3 \otimes y,\]
which yields  $2\pi\zeta(3)$.

If $y=w$,
\[\Omega=\{y\}_3\otimes y,\]
whose integral is zero.

When $y=1$, in this case, $w=1$ or $z=1$. If $w=1$, $\Delta=0$.

If $z=1$,
\[\Delta= \left\{ 1 - \frac{1}{x}\right\}_2 \otimes w \wedge (1-w) - \left\{1-\frac{1}{w} \right\}_2 \otimes \left(1-\frac{1}{x}\right) \wedge w - \{w\}_2\otimes x\wedge \left(1-\frac{1}{x}\right).\]
Only the term in the middle yields a nonzero differential form. In fact, the term in the middle yields
\[ \Omega= \{w\}_3\otimes \left( 1 - \frac{1}{x} \right).\]

$\Omega$ will be integrated on the boundary which is $x=1$, $w=1$ and $x=w$.

If $x=1$, $\Omega=0$.
If $w=1$,
\[\Omega = \{1\}_3 \otimes \left(1 - \frac{1}{x} \right).\]
This integration is  equal to $\pi \zeta(3)$.

If $x=w$,
\[ \Omega = \{x\}_3\otimes \left(1 - \frac{1}{x} \right),\]
which integrates to zero.

When $w=x$, (in this case, $z=y$ unless $x=1$).

\[\Delta  = -\{x\}_2 \otimes x \wedge y +  \left\{y-\frac{y}{x}\right\}_2 \otimes  x \wedge y(1-x)\]
\[ - \left\{y-\frac{y}{x}\right\}_2 \otimes \left(y-\frac{y}{x}\right) \wedge x
- \{x\}_2\otimes x\wedge \left(y-\frac{y}{x}\right). \]

Then
\[\Omega = -2\{x\}_3\otimes y-2 \left\{y-\frac{y}{x}\right\}_3 \otimes  x - \{x\}_3\otimes  \left(1-\frac{1}{x}\right). \]

Now $\Omega$ is to be integrated on the boundary, which is $x=1$, $y=1$ and $x=y$ (see picture \ref{cube}).

If $x=1$,
\[ \Omega = -2\{1\}_3\otimes y,\]
which gives $4\pi \zeta(3)$.

If $y=1$,
\[\Omega = -2 \left\{1-\frac{1}{x}\right\}_3 \otimes  x - \{x\}_3\otimes  \left(1-\frac{1}{x}\right). \]
Now use that
\[ \{x\}_3+\{1-x\}_3 +\left\{1- \frac{1}{x} \right \}_3 = \{1\}_3 \]
and the fact that $|x|=1$ to conclude
\[ -2 \left\{1-\frac{1}{x}\right\}_3 \otimes  x = \{x\}_3 \otimes x - \{1\}_3 \otimes x.\]
The total integration  in this case is $2 \pi \zeta(3)$.

If $x=y$,
\[ \Omega = -2\{x\}_3\otimes x-2 \left\{x-1\right\}_3 \otimes  x - \{x\}_3\otimes  \left(1-\frac{1}{x}\right)   \]
which leads to $-2 \oint \mu( x-1,x)$ (we will not need to compute this integral for the final result).

When $w=y$, (in this case, $z=x$ unless $y=1$),
\[ \Delta = -\{x\}_2 \otimes y \wedge x +  \left\{y-\frac{y}{x}\right\}_2 \otimes  y \wedge \left( x - \frac{x}{y}\right)\]
\[ -\left\{x-\frac{x}{y}\right\}_2 \otimes \left(y-\frac{y}{x}\right) \wedge y
- \{y\}_2\otimes x\wedge \left(y-\frac{y}{x}\right)   \]

\[= \{x\}_2 \otimes x \wedge y + \{y\}_2 \otimes y \wedge x - \{y\}_2\otimes x\wedge \left(1-\frac{1}{x}\right)\]
\[ -  \left\{y-\frac{y}{x}\right\}_2 \otimes  \left( x - \frac{x}{y}\right) \wedge y -\left\{x-\frac{x}{y}\right\}_2 \otimes \left(y-\frac{y}{x}\right) \wedge y.\]

By the five-term relation,
\[ \left\{ 1 - \frac{1}{x} \right\}_2 +\{y\}_2 + \left\{ 1 -y\left(1-\frac{1}{x}\right)\right\}_2 + \left\{\frac{1}{x+y-xy}\right\}_2 + \left\{ \frac{1-y}{1-y+\frac{y}{x}}\right\}_2 = 0\]
\[ \left\{ x \right\}_2 +\{y\}_2 - \left\{ y - \frac{y}{x}\right\}_2 - \left\{x+y-xy\right\}_2 - \left\{ x- \frac{x}{y}\right\}_2 = 0.\]

Then we obtain
\[ \Delta = \{x\}_2 \otimes x \wedge y + \{y\}_2 \otimes y \wedge x - \{y\}_2\otimes x\wedge \left(1-\frac{1}{x}\right)\]
\[ - \{x\}_2 \otimes  \left( x - \frac{x}{y}\right) \wedge y - \{y\}_2 \otimes  \left( x - \frac{x}{y}\right) \wedge y+ \{x+y-xy\}_2 \otimes  \left( x - \frac{x}{y}\right) \wedge y+\left\{x-\frac{x}{y}\right\}_2\otimes  \left( x - \frac{x}{y}\right) \wedge y\]
\[- \{x\}_2 \otimes \left(y-\frac{y}{x}\right) \wedge y - \{y\}_2 \otimes  \left( y - \frac{y}{x}\right) \wedge y+ \{x+y-xy\}_2 \otimes  \left( y - \frac{y}{x}\right) \wedge y+\left\{y-\frac{y}{x}\right\}_2\otimes  \left( y - \frac{y}{x}\right) \wedge y\]

\[ = \{x\}_2 \otimes x \wedge y + \{y\}_2 \otimes y \wedge x - \{y\}_2\otimes x\wedge \left(1-\frac{1}{x}\right)\]
\[ - \{x\}_2 \otimes  (1-x)(1-y)  \wedge y - \{y\}_2 \otimes  (1-x)(1-y) \wedge y +\left\{x-\frac{x}{y}\right\}_2\otimes  \left( x - \frac{x}{y}\right) \wedge y\]
\[ + \{x+y-xy\}_2 \otimes (1-x)(1-y ) \wedge y+\left\{y-\frac{y}{x}\right\}_2\otimes  \left( y - \frac{y}{x}\right) \wedge y.\]

Now
\[ - \{y\}_2\otimes x\wedge \left(1-\frac{1}{x}\right) - \{x\}_2 \otimes  (1-y)  \wedge y \]
is zero in the differential form.

Therefore,
\[\Omega =\{x\}_3 \otimes  y + \{y\}_3 \otimes x   +\{1-x\}_3 \otimes y   + \{y\}_3 \otimes  (1-x)(1-y) +\left\{x-\frac{x}{y}\right\}_3\otimes  y\]
\[- \{(1-x)(1-y)\}_3 \otimes  y+\left\{y-\frac{y}{x}\right\}_3\otimes  y.\]

$\Omega$ will be integrated on the boundary, which is $x=1$, $y=1$ and $x=y$.

If $x=1$,
\[\Omega= \{1\}_3 \otimes  y  +\{y\}_3\otimes(1-y)  +\left\{1-\frac{1}{y}\right\}_3\otimes  y\]
whose integral is  $3 \pi \zeta(3) $.

If $y=1$,
\[\Omega =  \{1\}_3 \otimes x +\{1\}_3\otimes(1-x) \]
which gives $3\pi \zeta(3)$.

If $x=y$,
\[\Omega =2 \{x\}_3 \otimes  x   +\{1-x\}_3 \otimes x   +2 \{x\}_3 \otimes  (1-x) +2 \left\{x-1 \right\}_3\otimes  x - \{(1-x)^2\}_3 \otimes  x\]
\[ = 2 \{x\}_3 \otimes  x  - 3\{1-x\}_3 \otimes x   +2 \{x\}_3 \otimes  (1-x) -2 \left\{x-1 \right\}_3\otimes  x,\]
which yields $3\pi \zeta(3) + 2 \oint \mu( x-1,x) $.

The poles are with $w=1$ but  $\Delta =0$ in this case. On the other hand, if $z=0$, then $w = x+y-xy$. But $|w|=1$ implies that $x=1$, $y=1$ or $x=-y$. In the first two cases, $w=1$ and $\Delta =0$. In the third case
\[ \Delta = -\{x\}_2 \otimes x^2 \wedge (1-x^2) -\{x^2\}_2 \otimes x \wedge (1-x),\]
which corresponds to zero if $|x|=1$.

Thus,
\[  8\pi^3m(P) = 36\pi  \zeta(3).\]

Finally,
\[m(P) = \frac{9}{2\pi^2} \zeta(3).\]
\end{proof}

\section{The $n$-variable case.}

The usual application of Jensen's formula (as in equation (\ref{eq:Jensen}) ) allows us to write, for $P \in \C[x_1, \dots, x_n]$,
\begin{equation}
m(P)  = m(P^*) + \frac{1}{(-2\pi\ii)^{n-1}} \int_G \eta_n(n)(x_1,\dots,x_n),
\end{equation}
where
\[G= \{P(x_1,\dots, x_n) =0\} \cap \{|x_1|=\dots=|x_{n-1}|=1, |x_n|\geq 1 \}.\]
(Recall that this is due to Deninger \cite{D}).

It is easy to see that we can then follow  a process that is analogous to the ones we followed for up to four variables. It remains, of course, to find an general algebraic way of describing the successive sets that we obtain by taking boundaries. Suppose that we do have a good description of the boundaries inside certain algebraic varieties, say $S_1 = \{P(x_1, \dots, x_n)=0\}, \dots, S_{n-1}$. Write, as usual, 
\[ \partial \gamma = \sum_k \epsilon_k [w_k], \quad \epsilon_k = \pm 1\]
where $\gamma$ is the collection of paths $S_{n-1} \cap \{|x_1|=1\}$. 
In principle, we should expect:
\begin{conj}
 Let $P(x_1, \dots, x_n) \in \R[x_1,...,x_n]$ be nonreciprocal. Assume that the following conditions are satisfied:
\begin{equation}\label{eq1}
x_1 \wedge \dots \wedge x_n = \sum_{i_1} r_{i_1} z_{i_1} \wedge (1-z_{i_1}) \wedge Y_{i_1}
\end{equation}
in $\bigwedge^n(\C(S_1)^*)\otimes \Q$, 
\begin{equation}\label{eq2}
\{z_{i_1}\}_2 \otimes Y_{i_1}  = \sum_{i_2} r_{i_1,i_2} \{z_{i_1,i_2}\}_2 \otimes z_{i_1,i_2} \wedge  Y_{i_1,i_2}
\end{equation}
in $(\mathcal{B}_2(\C(S_2)) \otimes \bigwedge^{n-2} \C(S_2)^*)_\Q$, for all $i_1$; More generally, assume that for $k=4, \dots, n-2$ we have 
\begin{equation}
\{z_{i_1,\dots i_{k-1}}\}_{k} \otimes Y_{i_1,\dots i_{k-1}}  = \sum_{i_k} r_{i_1,\dots, i_{k-1},i_k} \{z_{i_1,\dots, i_{k-1},i_k}\}_k \otimes z_{i_1,\dots, i_{k-1},i_k} \wedge  Y_{i_1,\dots,i_{k-1},i_k} 
\end{equation}
in $(\mathcal{B}_{k}(\C(S_{k})) \otimes \bigwedge^{n-k} \C(S_{k})^*)_\Q$, for all $i_1, \dots, i_{k-1}$. Finally, assume
\begin{equation}
\{z_{i_1, \dots, i_{n-2}}\}_{n-1} \otimes Y_{i_1, \dots, i_{n-2}}  = \sum_{i_{n-1}} r_{i_1, \dots, i_{n-2},i_{n-1}} \{z_{i_1, \dots, i_{n-2},i_{n-1}}\}_{n-1} \otimes z_{i_1, \dots, i_{n-2},i_{n-1}} 
\end{equation}
in $(\mathcal{B}_{n-1}(\C(S_{n-1})) \otimes \C(S_{n-1})^*)_\Q$, for all $i_1, \dots, i_{n-2}$.

Then we may write
\begin{equation}
 (2\pi)^{n-1}(m(P^*)-m(P)) =  \mathcal{L}_n(\xi)
 \end{equation}
 for
 \[ \xi = \sum_k \sum_{i_1, \dots,i_{n-1}} \epsilon_k r_{i_1} \dots r_{i_1, \dots,i_{n-1}} \{z_{i_1, \dots,i_{n-1}} (w_k)\}_n \]
\end{conj}

Here we have written $Y_{i_1}$, $Y_{i_1,i_2}$, ..., to denote elements in $\bigwedge^{n-2}(\C(S_1)^*)\otimes \Q$, $\bigwedge^{n-3}(\C(S_2)^*)\otimes \Q$, ...  A solution to equation (\ref{eq1}) determines the $Y_{i_1}$'s. Once the $Y_{i_1}$'s are defined, we solve equation (\ref{eq2}) and obtain the $Y_{i_1,i_2}$'s. The procedure continues in this fashion until we reach the $Y_{i_1,\dots,i_{n-2}}$'s.

Ideally, we would expect that this setting explains the nature of the $n$-variable examples described in \cite{L2}.

\subsection{The case of an $n$-variable family.}
Let us consider the case of the family of $n+1$-variable rational functions
 \[ z = \left( \frac{1-x_1}{1+x_1}\right) \dots \left( \frac{1-x_n}{1+x_n}\right),\]
whose Mahler measure was computed in \cite{L2}.

Though we are not able to perform all the steps for general $n$, we can at least prove that the first two differentials $\eta_{n+1}(n+1)$ and $\eta_{n+1}(n)$ are exact.

In this case the wedge product is
\[ x_1 \wedge \dots \wedge x_n \wedge z = \sum_{i=1}^n \left(x_1 \wedge \dots \wedge x_n \wedge (1-x_i) - x_1 \wedge \dots \wedge x_n \wedge (1+x_i) \right) \]
\[ = \sum_{i=1}^n (-1)^{i(n-1)} \left(x_i \wedge (1-x_i) \wedge x_{i+1}\wedge \dots \wedge x_{i+n-1} - x_i \wedge (1+x_i) \wedge x_{i+1}\wedge \dots \wedge x_{i+n-1} \right),\]
with the cyclical convention that $x_{i+n} = x_i$.

Thus we proved that $\eta=\eta_{n+1}(n+1)(x_1,\dots,x_n,z)$ is exact. The next step is to integrate $\eta_{n+1}(n)$ evaluated on the following element:
\[ \Delta =  \sum_{i=1}^n (-1)^{i(n-1)} ( \{x_i\}_2 \otimes x_{i+1} \wedge \dots \wedge x_{i+n-1} -
\{-x_i\}_2 \otimes x_{i+1} \wedge \dots \wedge x_{i+n-1}).\]

We are going to prove that $\omega=\eta_{n+1}(n)(\Delta)$ is exact. This form is defined in the variety $Z$ which is the projective closure of the algebraic set determined by
\[ (-1)^n = \left( \frac{1-x_1}{1+x_1}\right)^2 \dots \left( \frac{1-x_n}{1+x_n}\right) ^2.\]

We wish to show that $\omega$ is trivial in $H_{DR}^{n-1}(Z)$. First observe that 
\[ Z = Z_+ \cup Z_-,\]
where  $Z_{\pm}$ is given by the equation  
\[\pm z^*  = \left( \frac{1-x_1}{1+x_1}\right) \dots \left( \frac{1-x_n}{1+x_n}\right) \]
and
\[ z^* = \left\{ \begin{array}{cl} 1 & n \,\, \mbox{even}\\  \ii & n \, \, \mbox{odd} \end{array} \right. \]

In general, consider the variety given by the projective closure of the zeros of
\[\alpha=\left( \frac{1-x_1}{1+x_1}\right) \dots \left( \frac{1-x_n}{1+x_n}\right) \]
with $\alpha$ a nonzero complex number. This variety is birational to $\PP^{n-1}$ (this is easy to see by setting $y_i =  \frac{1-x_i}{1+x_i}$).

Hence we may think of each $Z_\pm$ as a copy of $\PP^{n-1}$. The singular points for this birational map are when $x_i=\pm 1$.

Now suppose that $n$ is even, $n=2k$.
 
If we prove that $\omega$ can be extended to the whole $Z_\pm$ and that this extension is consistent with the birationality of $Z_\pm$, it would imply that $\omega$ is closed, and that it could be seen as a class in $H_{DR}^{2k-1}(\PP^{2k-1})=0$ and then $\omega$ would be exact.
 
In order to extend $\omega$, we need to consider the points where some $x_i$ is equal to $1$, $-1$ (the points where the equation has singularities) and $0$, $\infty$  (the points where $\omega$ is not defined).

Consider the diagram
\[\begin{array}{ccc} {\mathcal B}_2(\C(Z)) \otimes \bigwedge^{n-1} \C(Z)^* & \stackrel{\eta_{n+1}(n)}{\longrightarrow} & \mathcal{A}^{n-1}(Z)(n)\\\\
\partial_v \downarrow &&\res_v \downarrow \\\\
{\mathcal B}_2(\C(Z)_v) \otimes \bigwedge^{n-2} \C(Z)_v^* & \stackrel{\eta_{n}(n-1)}{\longrightarrow} & \mathcal{A}^{n-2}(Z_v)(n-1)\end{array}\]
which describes the relation between the tame symbol and the residue morphism.

We would like to see that
\[ \res_v(\eta_{n+1}(n)(\{x_1\}_2\otimes x_2 \wedge \dots \wedge x_n)) =0,\]
where $v$ is the valuation defined by $x_i = \pm 1, 0, \infty$ for some $i$. Instead, we will see that 
\[\eta_n(n-1)(\partial_v(\{x_1\}_2\otimes x_2 \wedge \dots \wedge x_n))=0.\]

First suppose $x_i =1$. Then if $i \not = 1$, reducing modulo $x_i -1$ implies that $x_i =1$ and the only term that is possibly nonzero in $\partial_v(\{x_1\}_2\otimes x_2 \wedge \dots \wedge x_n)$ is $v(x_i) \{\bar{x}_1\}_2 \otimes \bar{x}_2\wedge \cdots \widehat{x}_i \dots \wedge \bar{x}_n$. However, $x_i$ is clearly not a uniformizer for $x_i -1$. Then the tame symbol is zero. If $i=1$, $\{x_1\}_2$ reduces to $\{1\}_2$ which corresponds to zero in $\eta_n(n-1)$, so we get zero again. The case $x_i = -1$ is analogous.

Now consider the case with $x_i=0$ for some $i$. Then it is easy to see that
\begin{equation}\label{eqtame}
\partial_v(\{x_1\}_2\otimes x_2 \wedge \dots \wedge x_n)  =   \left \{ \begin{array}{cc} \{\bar{x}_1\}_2 \otimes \bar{x}_2\wedge \cdots \widehat{x}_i \dots \wedge \bar{x}_n & \mbox{if}\quad i \not = 1 \\\\ 0 & \mbox{if}\quad i = 1 \end{array}\right .
\end{equation}
Since $x_i =0$, we are now in the variety defined by the projective closure of the zeros of the equation
\[  1 = \left( \frac{1-x_1}{1+x_1}\right)^2 \dots \widehat{\left( \frac{1-x_i}{1+x_i}\right)}^2 \dots  \left( \frac{1-x_n}{1+x_n}\right) ^2 .\]
We are in a situation that is analogous to the initial one. In other words, we are in the projective space $\PP^{2k-2}$. We would like to proceed by induction. In order to prove that  $\eta_{n-1}(n-2) (\{\bar{x}_1\}_2 \otimes \bar{x}_2\wedge \cdots \widehat{x}_i \dots \wedge \bar{x}_n)$ is trivial, we can prove that the tame symbols $\partial_w(\{\bar{x}_1\}_2 \otimes \bar{x}_2\wedge \cdots \widehat{x}_i \dots \wedge \bar{x}_n)$  are trivial by induction. However, $H_{DR}^{2k-2}(\PP^{2k-2})\cong \R$, so even if the symbols are trivial we will not be able to conclude that the form
$\eta_{n-1}(n-2) (\{\bar{x}_1\}_2 \otimes \bar{x}_2\wedge \cdots \widehat{x}_i \dots \wedge \bar{x}_n)$ is exact. What we can conclude is that it is either a generator for $H_{DR}^{2k-2}(\PP^{2k-2})$ or trivial. We would like to eliminate the first possibility.

Suppose, in order to make  notation easier, that $n=i$. Assume that $\eta_{n-1}(n-2) (\{\bar{x}_1\}_2 \otimes \bar{x}_2\wedge \dots \wedge \bar{x}_{n-1})$ is a generator for $H_{DR}^{2k-2}(\PP^{2k-2})$. By Poincar\'e duality, the integral
 \[ I=\int_{1=\left( \frac{1-x_1}{1+x_1}\right)^2 \dots \left(\frac{1-x_{2k-1}}{1+x_{2k-1}}\right)^2} \eta_{2k-1}(2k-2)( \{x_1\}_2 \otimes x_2\wedge \dots \wedge x_{2k-1}) \]
must be nonzero.

Now the transformation $x_i \rightarrow x_i^{-1}$ does not change the orientation of the variety but changes the sign of the differential $\omega$. Hence, $I= -I$ and that implies that $I=0$. Hence $\omega$ can not be a generator for $H_{DR}^{2k-2}(\PP^{2k-2})$ and it must be exact.

The case when $x_i = \infty $ is analogous.

Now suppose that $n=2k+1$ is odd. Then we may proceed as before. We have that $\omega$ can be seen as a class in $H_{DR}^{2k}(\PP^{2k})\cong \R$ and we can conclude that is exact by using the same idea that we used for the even case.

To conclude, $\omega=\eta_{n+1}(n)(\Delta)$ is exact and it must be the differential of certain $\mu =\eta_{n+1}(n-1)(\Omega)$. However, we were unable to find the precise formula for $\mu$. The results of \cite{L2} suggest that one should be able to continue this process to reach $\eta_{n+1}(1)$.

\bigskip
\begin{acknowledgements}

I would like to express my deepest gratitude to my Ph.D. supervisor, Fernando Rodriguez-Villegas, for his invaluable guidance and encouragement along this and other projects.  I am grateful to Vincent Maillot for many enlightening discussions about cohomology and the $n$-variable case.  I am also thankful to Herbert Gangl for helpful discussions about polylogarithms and Bloch groups. I wish to express my appreciation for the careful work of the referees, which has certainly improved the clarity and exposition of this paper.

This work is part of my Ph.D. dissertation at the Department of Mathematics at
the University of Texas at Austin. I am indebted to John Tate and the Harrington
fellowship for their generous financial support during my graduate studies.

This research was partially completed while I was a visitor at the Institut des Hautes \'Etudes Scientifiques. I am grateful for their support. 

\end{acknowledgements}


\begin{thebibliography}{BR-V03}

\bibitem[Be05]{Be} M. J. Bertin, Mesure de Mahler d'hypersurfaces $K3$, (preprint, January 2005).


\bibitem[Boy98]{B2} D. W. Boyd, Mahler's measure and special values of  L-functions, {\em Experiment. Math.\/} {\bf 7} (1998), 37--82.

\bibitem[BR-V02]{BRV1} D. W. Boyd, F. Rodriguez-Villegas, Mahler's measure and the dilogarithm (I), {\em Canad. J. Math.\/} {\bf 54} (2002), 468--492.

\bibitem[BR-V03]{BRV2} D. W. Boyd, F. Rodriguez-Villegas, with an appendix by N. M. Dunfield, Mahler's measure and the dilogarithm (II), (preprint, July 2003).



\bibitem[DL06]{DL} C. D'Andrea, M. N. Lal\'{\i}n,  On The Mahler measure of resultants in small dimensions, (to appear in {\em J. Pure Appl. Algebra}).

\bibitem[Den97]{D} C. Deninger,  Deligne periods of mixed motives, $K$-theory and the entropy of certain $ Z\sp n$-actions, {\em  J. Amer. Math. Soc. \/} {\bf 10}  (1997),  no. 2, 259--281.




\bibitem[Gon95a]{G3} A. B. Goncharov, Geometry of Configurations, Polylogarithms, and Motivic Cohomology, {\em Adv. Math.\/}  {\bf 114}  (1995),  no. 2, 197--318.



\bibitem[Gon02]{G4} A. B. Goncharov, Explicit regulator maps on polylogarithmic motivic complexes.  {\em Motives, polylogarithms and Hodge theory, Part I (Irvine, CA, 1998),\/}  245--276, {\em Int. Press Lect. Ser., 3, I, Int. Press, Somerville, MA,} 2002.


\bibitem[Gon05]{G5} A. B. Goncharov, Regulators. {\em Handbook of $K$-theory}. Vol. 1, 2,  295--349, Springer, Berlin, 2005.






\bibitem[Lal06]{L2} M. N. Lal\'{\i}n, Mahler measure of some n-variable polynomial families, {\em J. Number Theory\/}  {\bf 116} (2006), no. 1, 102--139.

\bibitem[Lal07]{L3} M. N. Lal\'{\i}n, Mahler measures and computations with regulators (in preparation).



\bibitem[R-V99]{RV} F. Rodriguez-Villegas, Modular Mahler measures I, {\em Topics in number theory (University Park, PA 1997),\/} 17--48, {\em Math. Appl., 467, Kluwer Acad. Publ. Dordrecht,\/} 1999.



\bibitem[Smy81]{S1} C. J. Smyth, On measures of polynomials in several  variables, {\em Bull. Austral. Math. Soc. Ser. A\/} {\bf 23} (1981), 49--63. Corrigendum (with G. Myerson): {\em  Bull. Austral. Math. Soc.\/} {\bf 26} (1982), 317--319.










\bibitem[Zag90]{Z1} D. Zagier, The Bloch-Wigner-Ramakrishnan polylogarithm function. {\em Math. Ann.\/} {\bf 286}  (1990),  no. 1--3, 613--624.

\bibitem[Zag91a]{Z0} D. Zagier, Special values and functional equations of polylogarithms.  Structural properties of polylogarithms,  377--400, {\em Math. Surveys Monogr., 37\/}, Amer. Math. Soc., Providence, RI, 1991

\bibitem[Zag91b]{Z2} D. Zagier, Polylogarithms, Dedekind Zeta functions, and the Algebraic $K$-theory of Fields,  {\em Arithmetic algebraic geometry\/} (Texel, 1989),  391--430, Progr. Math., 89, Birkh\"auser Boston, Boston, MA, 1991.

\bibitem[ZG00]{GZ} D. Zagier, H. Gangl, Classical and elliptic polylogarithms and special values of $\Lf$-series, {\em The arithmetic and geometry of algebraic cycles\/} (Banff, AB, 1998), 561 - 615, NATO Sci. Ser. C Math. Phys. Sci., 548, Kluwer Acad. Publ., Dordrecht, 2000.



\end{thebibliography}
\end{document}